\documentclass[10pt,letterpaper]{article}
\usepackage{amsmath}
\usepackage{caption}
\usepackage{subcaption}
\usepackage{xcolor}
\usepackage{dirtytalk}
\usepackage{times}
\usepackage{draftwatermark}
\SetWatermarkText{Draft}
\SetWatermarkScale{1.5}

\begin{document}
\thispagestyle{empty}
\pagestyle{empty}
 \vspace{\baselineskip}
 \vspace{\baselineskip}
\begin{center}
	\Large{\textbf{On the use of Nonlinear Normal Modes for Nonlinear Reduced Order Modelling}}\\
	\normalsize
	\vspace{\baselineskip}
	\vspace{\baselineskip}
	\vspace{\baselineskip}
	\vspace{\baselineskip}
	\large Thomas Simpson$^{1}$, Nikolaos Dervilis$^{2}$ and Eleni Chatzi$^{1}$\\ 
	\vspace{\baselineskip}
	$^{1}$ Institute of Structural Engineering, Department of  Civil, Environmental and Geomatic Engineering, ETH Zürich, Stefano-Franscini Platz 5, 8093, Zürich, Switzerland\\
	$^{2}$ Dynamics Research Group, Department of Mechanical Engineering, University of Sheffield.\\Sheffield S1 3JD, UK.\\
\end{center}
\normalsize
\vspace{\baselineskip}
\vspace{\baselineskip}
\vspace{\baselineskip}

\section*{ABSTRACT}
\label{abstract}
In many areas of engineering, nonlinear numerical analysis is playing an increasingly important role in supporting the design and monitoring of structures. Whilst increasing computer resources have made such formerly prohibitive analyses possible, certain use cases such as uncertainty quantification and real time high-precision simulation remain computationally challenging. This motivates the development of reduced order modelling methods, which can reduce the computational toll of simulations relying on mechanistic principles. The majority of existing reduced order modelling techniques involve projection onto linear bases. Such methods are well established for linear systems but when considering nonlinear systems their application becomes more difficult. Targeted schemes for nonlinear systems are available, which involve the use of multiple linear reduction bases or the enrichment of traditional bases. These methods are however generally limited to weakly nonlinear systems. In this work, nonlinear normal modes (NNMs) are demonstrated as a possible invertible reduction basis for nonlinear systems. The extraction of NNMs from output only data using machine learning methods is demonstrated and a novel NNM-based reduced order modelling scheme introduced. The method is demonstrated on a simulated example of a nonlinear 20 degree-of-freedom (DOF) system.

\section{INTRODUCTION}

The utilisation of nonlinear finite element simulations for dynamic problems has become more common in recent years as a result of increasing computational power rendering possible analyses that were previously prohibitively expensive. However, increasing interest in uncertainty quantification and digital twins once again raises the issue of computational resources. As such, the development of efficient metamodelling schemes for nonlinear dynamical systems is necessary. Meta-models are models which approximate a full-order or high-resolution model to some degree of accuracy, whilst greatly decreasing the cost of evaluation.

Reduction of FE models was initially considered with regards to the reduction of
linear elastic models in the 1960s and techniques such as the Craig-Bampton and
MacNeal-Rubin methods \cite{Craig1968}\cite{DeKlerk2008}, are well developed and
widely used for the reduction of linear FE models. These methods largely make use of
linear normal modes as a reduced basis on which the systems' equations of motion are projected. Normal modes form a very efficient reduction basis for linear systems due to their orthogonality and the ability to target certain frequency regions of interest. When considering nonlinear systems however, these linear methods are not appropriate, as the concept of linear normal modes no longer holds. The most common method for the reduction of weakly nonlinear FE models adopts Proper Orthogonal Decomposition (POD), wherein a dynamic simulation of the system is carried out with “snapshots” of the response field extracted at various time intervals. These snapshots are used to construct an optimal linear reduction basis using singular value decomposition \cite{Carlberg2009}. Related recent work explores the extension of the traditional reduction approach by inclusion of higher order enrichments, such as mode shape derivatives \cite{Wu2016}.

An alternative to this reduction approach is offered by a theory which extends the concept of normal modes to nonlinear systems, the so-called nonlinear normal modes (NNMs). Several different formulations of NNMs exist with the Shaw-Pierre modes being considered in this work, which build on the fundamental work of Rosenberg \cite{Rosenberg1962}. As a result of their ability to form a nonlinear manifold, these NNMs can provide a very efficient reduction basis for nonlinear systems \cite{Haller2016}. Furthermore, they should be generally valid over the whole phase space of a system and not require multiple, energy dependent reduction bases as in the case of linear methods. In the following work, NNMs are used as basis to construct a reduced order model of a nonlinear system. Autoencoder neural networks are used to extract these NNMs, a recurrent neural network (RNN) is then used to create a regression approximating the overall system, making use of the extracted NNMs as a basis.

The layout of this paper is as follows. Section 2 introduces the theoretical background of nonlinear normal modes along with the process of their extraction. Section 3 describes autoencoder neural networks, their architecture and common usage. Section 4 then discusses the conceptual framework of the reduced order modelling procedure used herein. Section 5 discusses the regression problem in the reduced order model and briefly introduces the long short term memory (LSTM) neural network used. Section 6 gives a demonstration of the reduced order modelling method in practice, the nonlinear MDOF system approximated is described along with greater details on the implementation of the ROM on this system. Finally section 7 contains the conclusions and further work considerations.

\section{Nonlinear Normal Modes}

\subsection{Shaw-Pierre Manifolds}
Shaw and Pierre first proposed their theory of non-linear normal modes in 1993. They defined a NNM as being a "a motion which takes place on a two dimensional invariant manifold in phase space". More simply put, this means that in a Shaw-Pierre NNM the motion of each point in the system could be given as a function of the displacement and velocity of a single point in the system. This function if plotted in phase space, would form a surface which is the manifold they refer to. In the case of the linear system, the manifold would be a flat plane. Nonlinearities in a system cause curvature of this manifold. This curvature results in the amplitude dependence of mode shapes observed in nonlinear systems \cite{Shaw1993}.

The calculation of these NNMs however is rather challenging, the assumption of the normal mode forming an invariant manifold is represented in equation \ref{eq:SP1}. In this equation $x_i$ and $\dot{x_i}$ are the displacement and velocity of the ith coordinate in the response whilst $u_i$ and $v_i$ represent the displacement and velocity of the ith transformed coordinate in the NNM. It can be seen that velocity and displacement of each of the coordinates in the NNM's can be described as a function of a displacement velocity pair of a single coordinate. These are linked by the functions $f_i$ and $g_i$. This relationship is then substituted into the equation of motion which results in a new formulation of the equations of motion, which are unfortunately at least as difficult to solve as the original case. The solution here is to assume a power series solution for the functions $f_1\cdots f_n$ and $g_1\cdots g_n$. This then allows for analytical solution for some cases and numerical solutions for the general case.

\begin{equation}\label{eq:SP1}
\begin{bmatrix}
x_1\\\dot{x}_1\\x_2\\\dot{x}_2\\\vdots\\x_n\\\dot{x}_n
\end{bmatrix}=
\begin{bmatrix}
u_1\\v_1\\u_2\\v_2\\\vdots\\u_n\\v_n
\end{bmatrix}=
\begin{bmatrix}
u_1\\v_1\\f_2(u_1,v_1)\\g_2(u_1,v_1)\\\vdots\\f_n(u_1,v_1)\\g_n(u_1,v_1)
\end{bmatrix}
\end{equation}

\subsection{Output Only Extraction}
Some previous work examined the use of machine learning methods to extract NNMs from output only data of nonlinear systems in a manner similar to how singular value decomposition can be used for principal orthogonal decomposition. The work of Worden and Green \cite{Worden2016} made use of the representation of NNMs as truncated polynomial series as defined in the Shaw-Pierre paper. Having assumed the form of these polynomials, the coefficients were inferred through an optimisation procedure wherein statistical independence up to a given order of correlation was also maximised. In addition, in the work of Dervilis et al, \cite{Dervilis2019}, various machine learning methods drawn from the field of manifold learning were demonstrated for the extraction of NNMs. Both these methods were illustrated on systems of few degrees of freedom for both simulated datasets and an experimental 4 DOF nonlinear dataset from Los Alamos national laboratory \cite{Figueiredo2009}.

\section{Autoencoder Neural Networks}
\begin{figure}[h!]
    \centering
    \includegraphics[width=100mm]{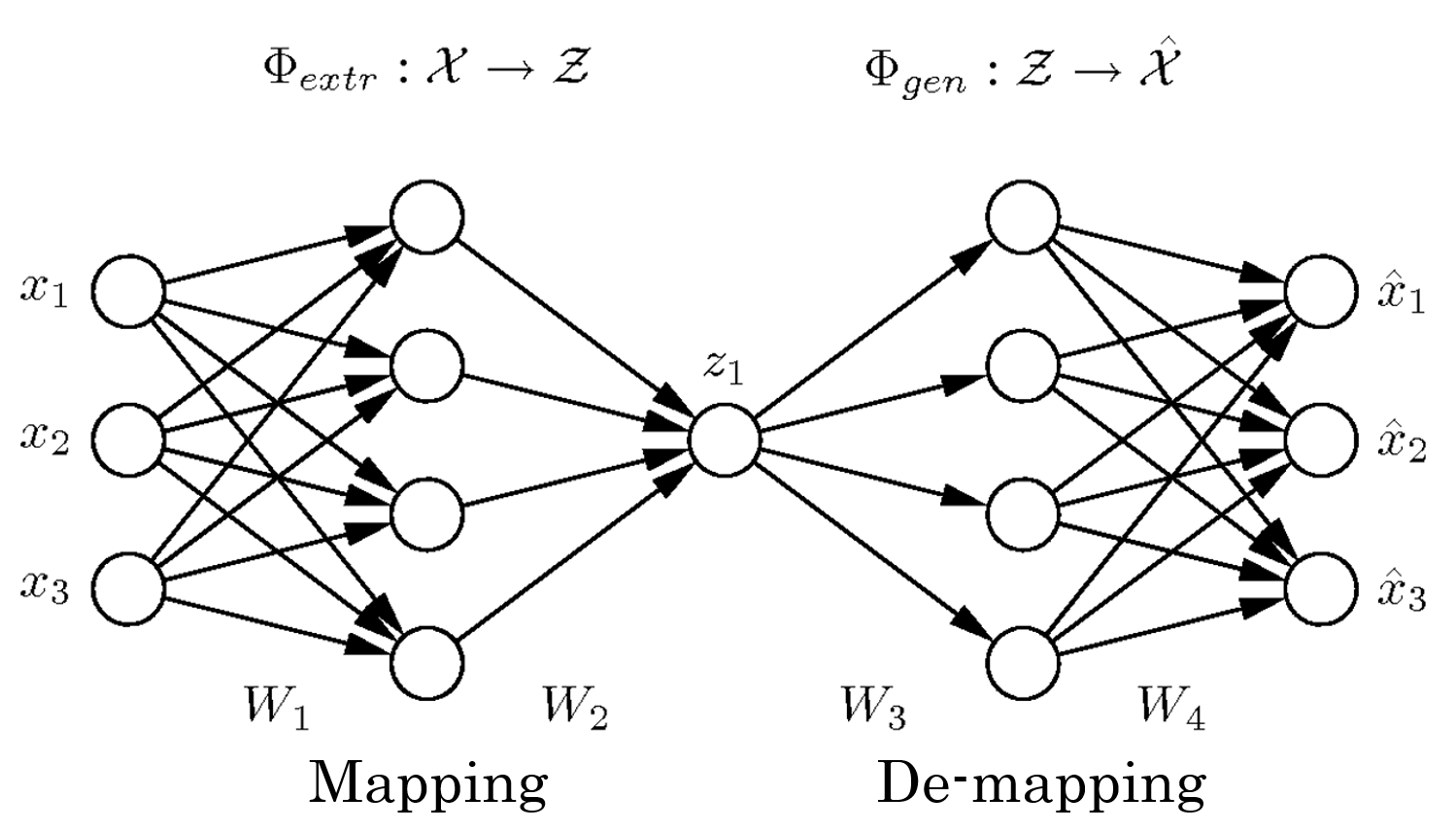}
    \caption{Architecture of an AE with 3 dimensional input/output and a 2 dimensional bottleneck.}
    \label{fig:AE}
\end{figure}

Autoencoder (AE) neural networks are a specialised type of neural network often used for dimensionality reduction or de-noising problems \cite{Hinton2006}. They are constructed, as illustrated in Figure \ref{fig:AE}, as a deep neural network architecture with a so-called "bottleneck" layer with a reduced number of nodes; in the case of this example to one. The key concept of their operation involves setting the cost function of the network such that the output attempts to recreate the inputs as closely as possible. By the inclusion of the bottleneck layer which forces the data through a lower dimensional feature space, a near optimal reduction of the data onto the feature space dimension can be achieved. The hidden layers in the network allow for a nonlinear transform of the input inputs before the bottleneck layer. The nature of an AE is such that it simultaneously learns both the encoding operation to transform from the physical to latent space, and the decoding operation whereby the latent variables are returned to the physical space.

Autoencoders have previously been used in similar applications in which it is desired that a model for a system be constructed on a lower dimensional "manifold", whereby it is believed that high dimensional data can be well described as existing on a lower dimensional space \cite{Yoo2017}. In the case of nonlinear structural dynamics it is thought that the NNMs could correspond to this lower dimensional manifold on which the majority of the high dimensional data can be represented. The same way that a high dimensional linear finite element model can often be represented on a far lower number of mode shapes. The encoding of an autoencoder is the function through which the high dimensional input data can be transformed to the lower dimensional manifold or NNM space, whilst the decoder can provide the inverse operation.

\section{Model Order Reduction with Autoencoders}
\begin{figure}[h!]
\centering
\includegraphics[width=160mm]{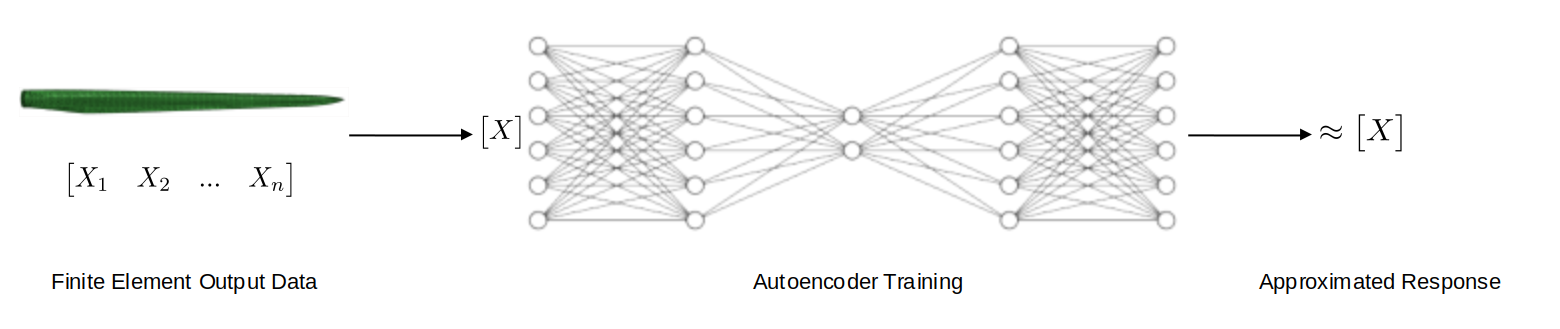}
\caption{1st step of ROM; a reduced number of NNMs are extracted from output only dataset}
\label{fig:ROM1}
\end{figure}

The first step of the model order reduction method used herein involves generating training data from the system of interest. This involves creating simulation time histories for the model of both the forcing and response of the system at all DOFs. Having generated time histories, the autoencoder can be used to compress the dimensionality down to a reduced number of latent variables or NNMs. These are extracted from the displacement time histories of the data with a number of NNMs retained, which balances reduction in the dimensionality of the system with accuracy of reconstruction of the full field response, as shown in Figure \ref{fig:ROM1}. In this process the cost function used is described by Equation \ref{eq:MSE}, which can be seen as the mean squared error of reconstruction between the original time series and those after compression.

\begin{equation}\label{eq:MSE}
    \ell(\hat{x})=\frac{1}{N}\sum_{i=1}^{N}(x_i-\hat{x_i})^2
\end{equation}

\begin{figure}[h!]
\centering
\includegraphics[width=140mm]{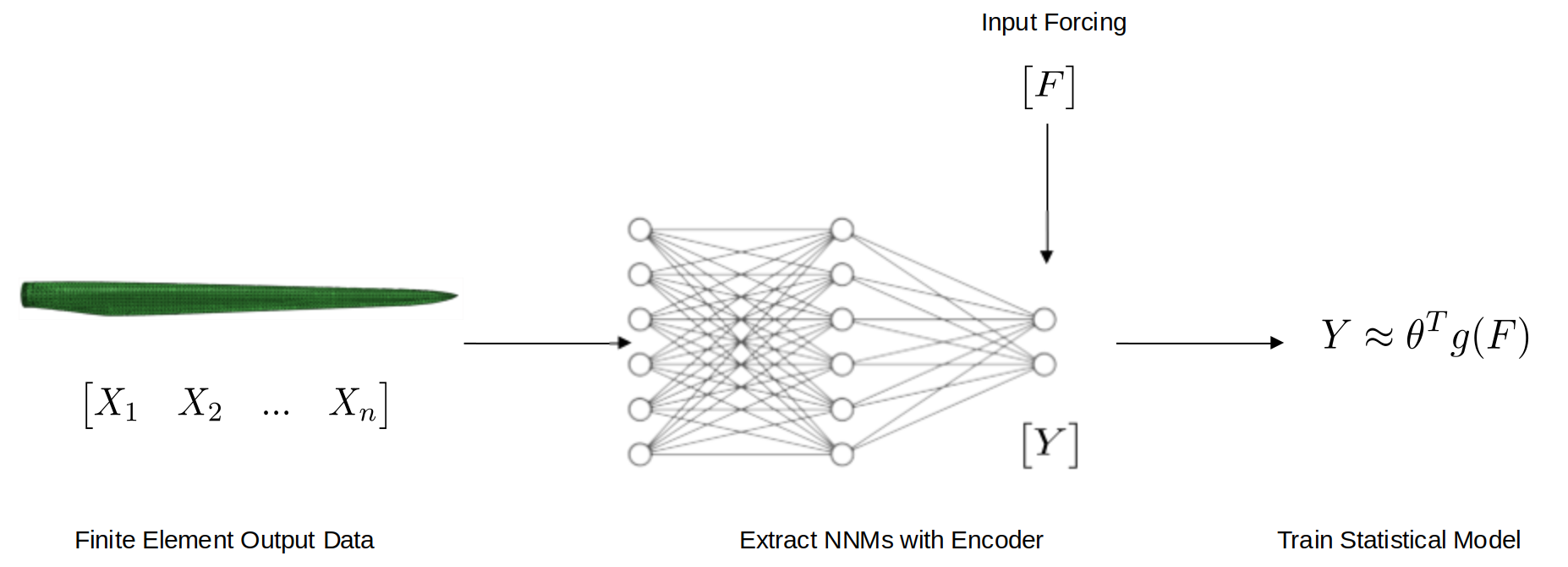}
\caption{2nd step of the ROM; a regression model is trained which predicts response in the retained NNMs based on forcing histories.}
\label{fig:ROM2}
\end{figure}

As opposed to projection based reduction methods, whereby the original equations of motion can be projected on to the reduced coordinate set and solved through conventional means, the method developed herein uses a statistical regression to allow for future simulations of the system. The regression model is trained which predicts the response of the system within the NNMs based on a given forcing time history, as illustrated in Figure \ref{fig:ROM2}. The regression model is trained based on the forcing time histories used for the simulation, $F$, and the corresponding NNMs are extracted by the encoder, $Y$.

\begin{figure}[h!]
\centering
\includegraphics[width=140mm]{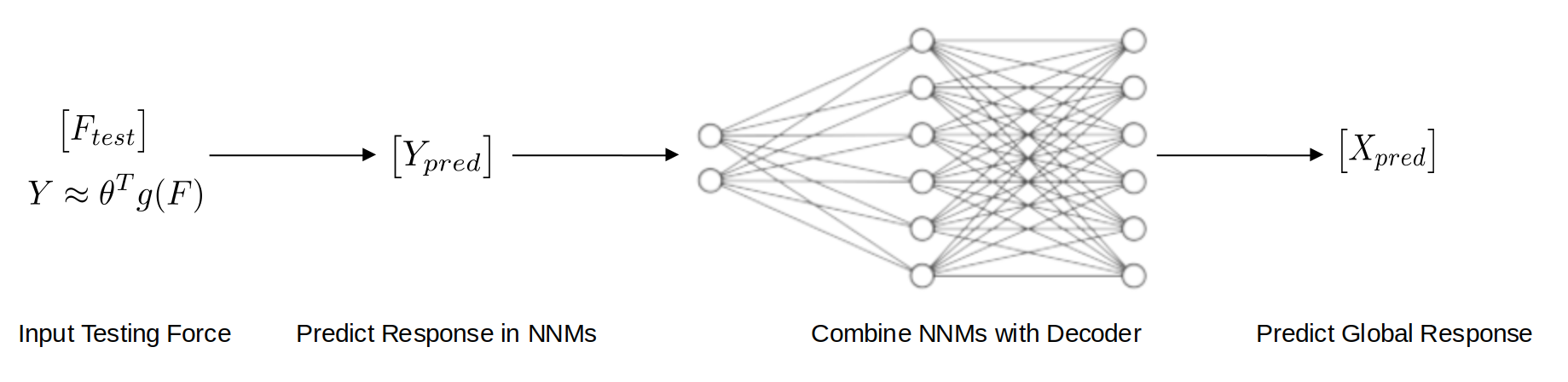}
\caption{3rd step of the ROM; the response in the NNMs is predicted for new forcing histories and the full response is reconstructed through the decoder}
\label{fig:ROM3}
\end{figure}

The final step in the metamodelling process is to predict the system response to a new forcing time history. As illustrated in Figure \ref{fig:ROM3}, the trained regression model is used to predict the response of the system within the NNMs, having predicted this response, $Y_{pred}$, the full field response, $X_{pred}$, can be recovered by using the decoding portion of the autoencoder.

\section{LSTM Neural Networks}
In constructing the meta-model, a statistical regression model must be used which can predict new time histories of the latent variables based on forcing time histories. The time series nature of the data would encourage the use of a method including an auto-regressive element. The nonlinear nature of the system would further require a nonlinear regression method to accurately represent the nonlinear dynamics. In this case it was decided to use a neural network regression architecture, the advantages of neural networks being the scalability and ability to perform MIMO regression.

Recurrent neural networks (RNNs) are a class of networks developed specifically for dealing with sequence data such as time series \cite{Ruineihart1985}. The core of an RNN is not dissimilar from that of an ARX model. A basic RNN consists of a unit which receives as input exogenous variables along with a hidden state value which is propagated from the previous time step. Based on these inputs a new value of the hidden state is calculated using a conventional neural network layer, i.e. through multiplication with a weight matrix and then application of an activation function. The output at each time step can then be calculated through an additional layer or layers which act upon the hidden state of the RNN. The hidden state of the RNN can preserve some context information about the previous states and inputs. In contrast to an ARX method however, it is not necessary to state apriori, how many previous values of input and outputs should be considered. This means that theoretically, RNNs can be very powerful at learning relationships between time series with long term dependencies, however, in practice it was found that when training an RNN using back propagation through time, the process suffered from vanishing gradients. As a result of this, a conventional RNN struggles to propagate errors back through many time steps and hence fails to learn long term dependencies \cite{Bengio1994}. This motivated the development of more advanced RNN cells, such as the long short term (LSTM) cell \cite{Hochreiter1997}.

Much of the significant work with RNNs has been achieved making use of LSTM cells. These have found very wide use in the field of natural language processing (NLP) \cite{Sundermeyer} as well as time series forecasting \cite{Gensler2017}. In the area of mechanical systems however, the use of LSTM networks is less prevalent with a limited amount of work applying them to traditional system identification problems \cite{Wang2017}. They have proven themselves however, to be capable of learning long term nonlinear relationships in NLP tasks and hence have significant potential as a useful regression framework for nonlinear time series.

\section{Case Study}
\begin{figure}[h!]
\centering
\includegraphics[width=150mm]{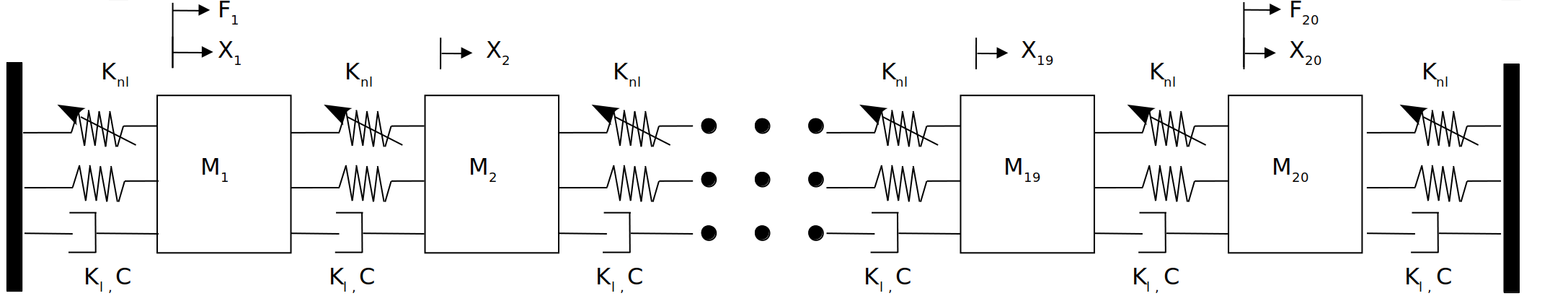}
\caption{The simulated 20 DOF system including cubic nonlinear spring elements between each mass.}
\label{fig:20dof}
\end{figure}

We verify applicability of the proposed method on a 20 degree-of-freedom mass, spring, damping system. The system was considered to be comprised of 20 masses in a chain, with linear spring, linear damping and nonlinear spring elements between each mass, as demonstrated in Figure \ref{fig:20dof}. The nonlinearity of the non-linear spring elements is of the cubic type. Forcing was applied to this system in two locations, at the first and final masses in the chain. All the parameters are assumed uniform across all degrees of freedom and the underlying linear system is proportionally damped; the parameters are defined as follows:  $M_i = 0.1\,kg$, $K_l = 100\,N/m$, $C_l = 0.1\,kg/s$, $K_{nl} = 2500\,N/m^{3}$. 

\begin{figure}[h!]
    \centering
    \includegraphics[width=150mm]{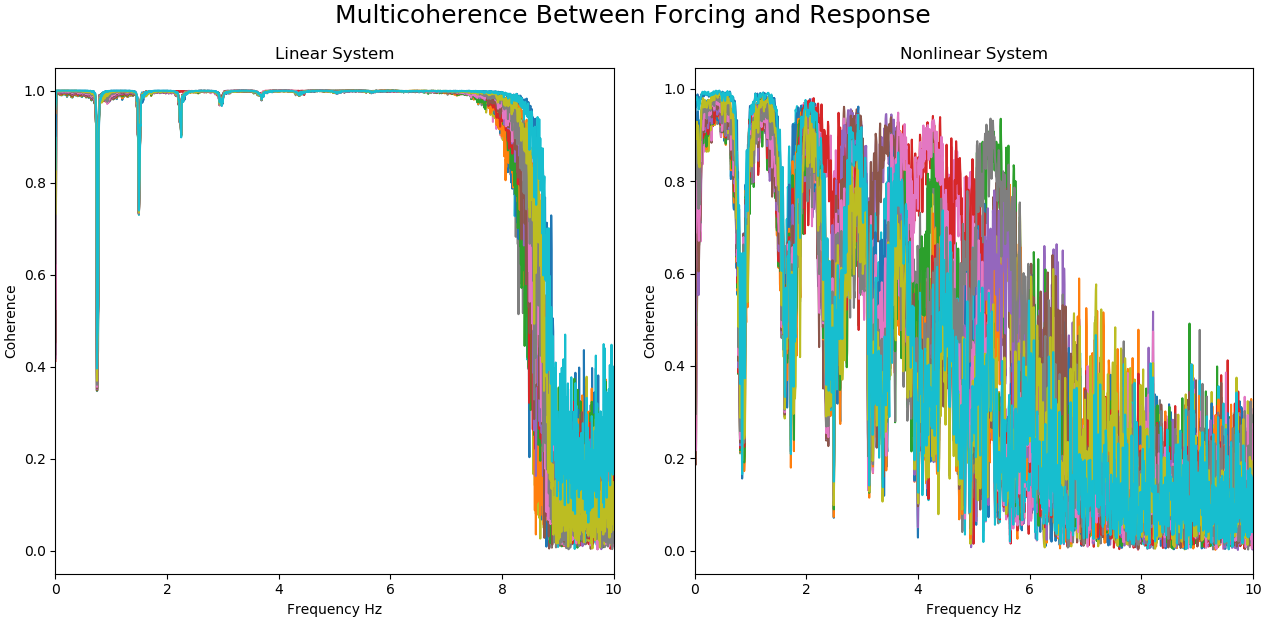}
    \caption{Multicoherence of the signal}
    \label{fig:Multicohere}
\end{figure}

Time series simulations of this system were carried out using a 4th order Runge-Kutta fixed step integration scheme. A sampling frequency of 100 Hz was used for the simulation with a total time series length of 1000s.  The forcing applied at masses 1 and 20 comprises realisations of a constant amplitude Gaussian white noise input, which had been low pass filtered with a cut-off frequency of 8 Hz. Figure \ref{fig:Multicohere} shows the multicoherence between the input forcing and the outputs of the system both for the underlying linear system and the nonlinear system tested. The nonlinear system shows considerably more incoherence with the input than the linear system, this demonstrates the considerable nonlinearity of the response. For both systems the region of interest is principally at frequencies below 8 Hz. In both cases the response begins to become completely incoherent at higher frequencies it is thought that this is largely an artifact of the processing, since the forcing contains very little energy at higher frequencies both the input and output in this region is very small making the coherence value very sensitive to numerical and processing errors.

\begin{figure}[h!]
    \centering
    \includegraphics[width=150mm]{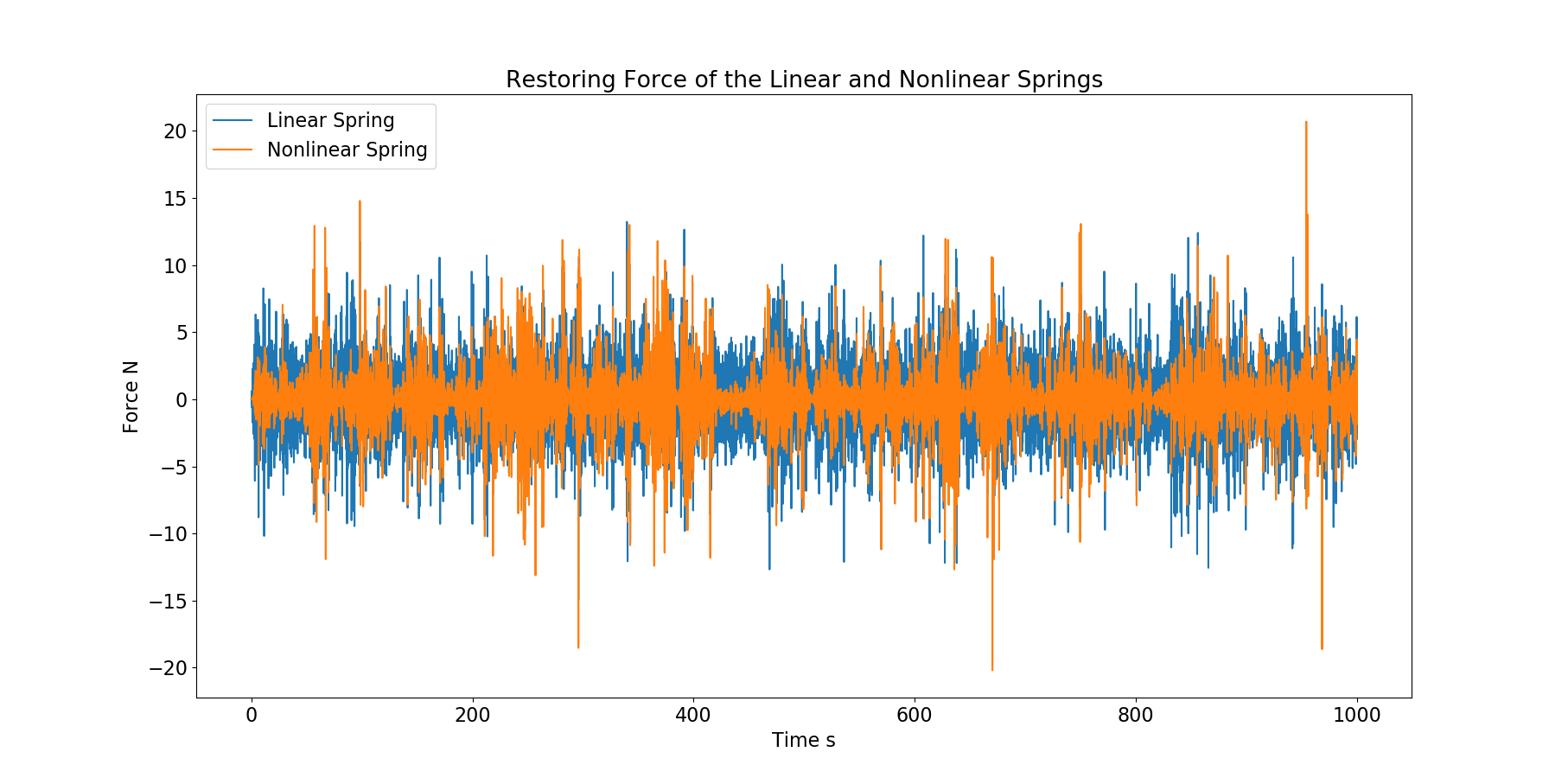}
    \caption{Restoring Force time History of the Linear and Nonlinear Spring Elements at the $10^{th}$ DOF}
    \label{fig:my_label}
\end{figure}

Figure \ref{fig:my_label} demonstrates the comparative magnitude of the restoring force in the linear and nonlinear spring elements of the system a the $10^{th}$ DOF. The magnitude of the nonlinear restoring force can be seen to be considerable regularly exceeding that of the linear elements. As such it can be considered that the system considered differs significantly from the underlying linear system. The remaining DOFs also demonstrated similar behaviour.

\subsection{NNM Extraction}
A 5 layer autoencoder is used for the extraction of NNMs from the output data. The encoder receives as input the 20 dimensional data corresponding to the displacement time histories at each of the 20 masses. The encoder portion consists of a hidden layer containing 20 nodes with a linear activation function, followed by a second 20-node layer with a hyperbolic tangent activation function. These layers were followed by the bottleneck layer consisting of 10 nodes. The decoder portion of the autoencoder mirrors the encoder in reverse, passing through a 20 node hidden layer with a hyperbolic tangent activation to a 20 node layer with a linear activation to the 20 dimensional output layer.

\begin{figure}[h!]
    \centering
    \includegraphics[width=140mm]{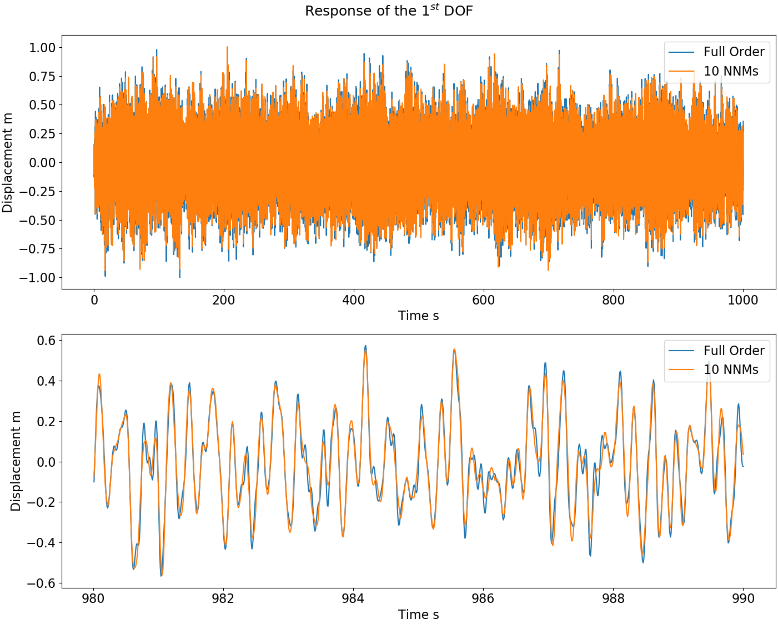}
    \caption{Reconstruction of the response at the 1st DOF with 10 retained NNMs}
    \label{fig:Reconstruct}
\end{figure}

For the training of the autoencoder, a training dataset consisting of the first 50 \% of the time histories was used corresponding to the first 500 s. The second 50 \% of the time histories were used as a testing dataset. It is noteworthy that the resulting training and testing loss values were similar, indicating that the NNMs captured were generalisable over the whole dataset. Figure \ref{fig:Reconstruct} shows the full order response along with the response reconstructed with 10 retained NNMs for the 1st DOF in the model.

\begin{figure}[h!]
    \centering
    \includegraphics[width=100mm]{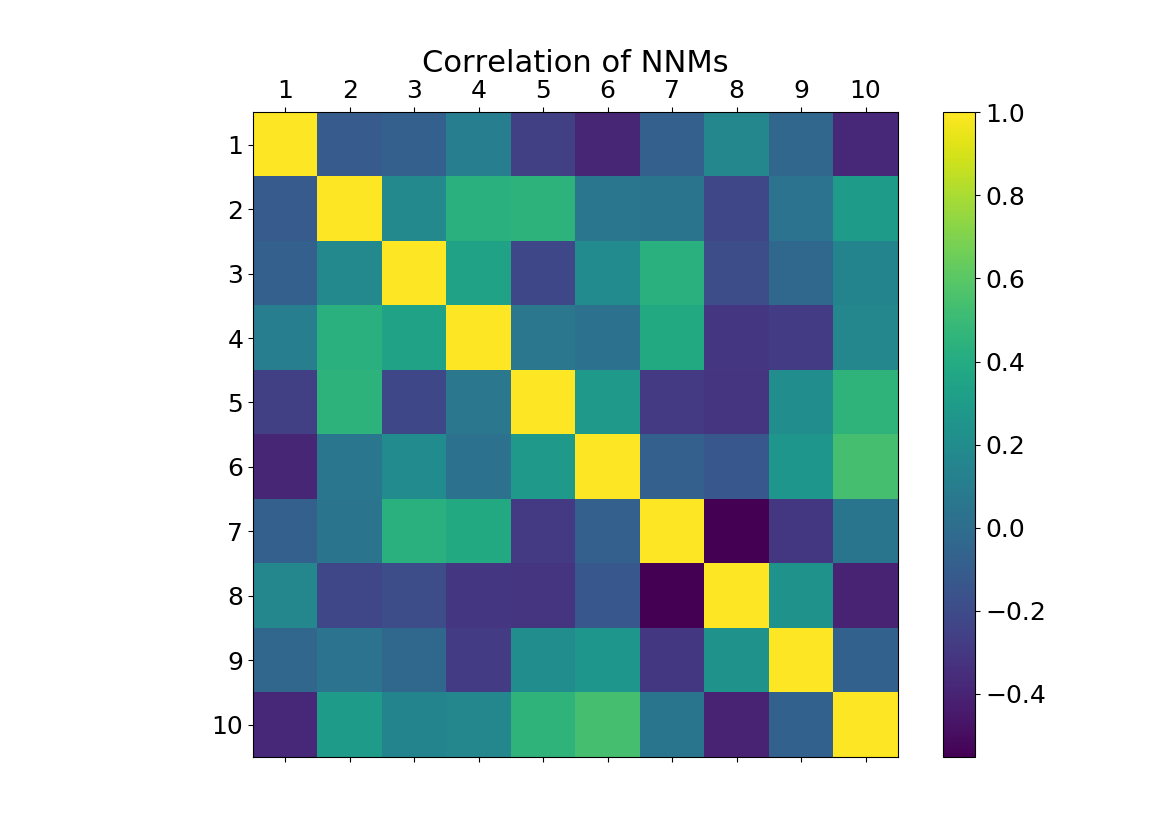}
    \caption{Correlation between the 10 extracted NNMs}
    \label{fig:NNMcorr}
\end{figure}

It is worth noting, that in previous studies involving output only extraction of NNMs using machine learning, that a key justification in referring to the extracted components as NNMs, is that the components be statistically independent. This statistical independence is considered to be analogous to the modal invariance property of theoretical NNMs \cite{Worden2016}\cite{Shaw1993} and the orthogonality property of linear modes. Figure \ref{fig:NNMcorr} presents the correlation matrix of the NNMs extracted by the autoencoder used herein. The non-zero off diagonal elements of this matrix imply significant correlation between the extracted components which implies that these are rather NNM-like quantities and do not comply with the strict definition of NNMs.

\subsection{Response Prediction in the Latent Space}
The response prediction in the NNMs involves a non-linear multivariate time series forecasting problem. This problem proved to be the most challenging aspect of the ROM process. Various regression architectures were attempted including polynomial NARX, Gaussian process NARX and recurrent neural net-works (RNN). Finally, a long short term (LSTM) neural network \cite{Hochreiter1997}, a type of RNN with improved long term memory characteristics, was used for this process. The regression was setup as an autoregressive process. This implies that the input for the prediction model at each time step includes not only the current and previous values of the forcing, but also previous values of the output. In training a prediction error was used, meaning that for each time step the ground truth values of the previous system outputs were given. During the testing stage however, a simulation error was used, meaning that the outputs predicted by the model at previous steps were subsequently fed to the model as inputs at a next time instant. The LSTM network used 64 hidden cells with 100 previous input and output values also fed as inputs to the network. The training data was assumed as the first 60\% of the time histories with predictions made on the next 1000 points corresponding to 10s, the exogenous inputs for this prediction were taken as the next 1000 points of forcing time history, which were not used in the training process. Initial states used for the first model prediction, were taken from the final “true” values in the training dataset.
\begin{figure}[h!]
    \centering
    \includegraphics[width=140mm]{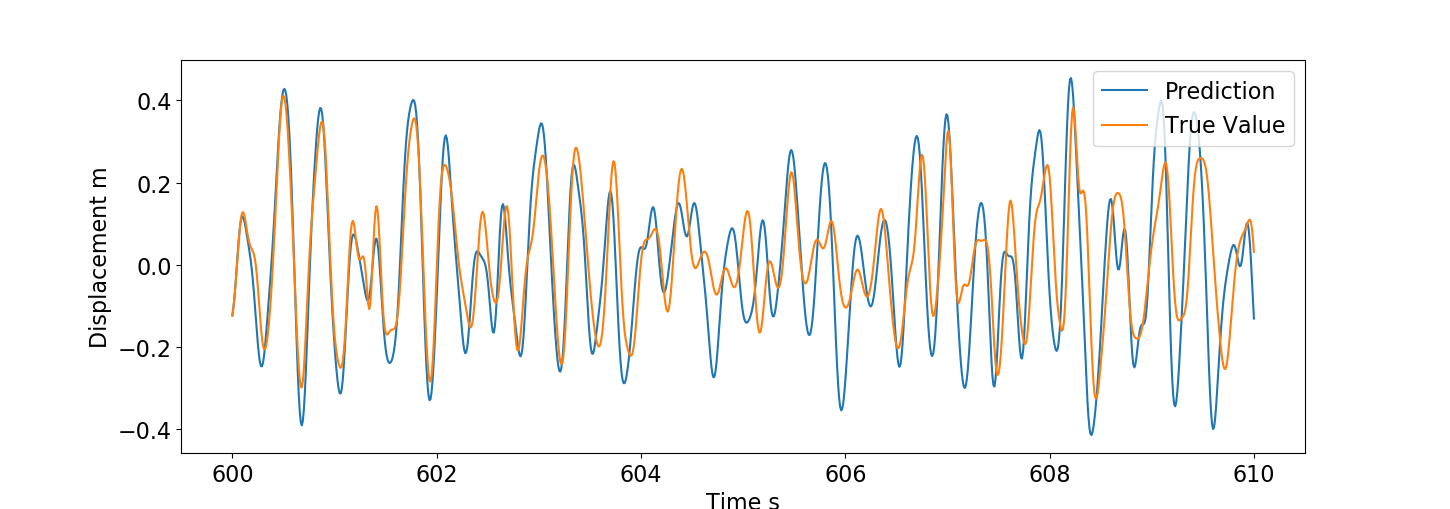}
    \caption{The respone prediction in the latent space for the 1st NNM.}
    \label{fig:NNMpredict}
\end{figure}

Figure \ref{fig:NNMpredict} compares the prediction of the response in one of the NNMs using the LSTM regression to the “true value”. The latter is the value obtained when the physical coordinate time series is passed through the encoder section of the autoencoder. The prediction performs well initially but diverges somewhat over time, it is noteworthy however that the prediction does not become unstable. Similar performance was obtained for all of the 10 NNMs predicted.

\subsection{Response Prediction in the Physical Space}

Upon availability of the response predictions of the system in the NNM variables, the full response prediction can be determined by simply passing the predicted NNM time histories through the decoder portion of the autoencoder. The decoder combines the NNMs so as to recreate the response in the physical space and is hence analogous to the nonlinear super-position discussed by Shaw and Pierre.

Figure 7 compares the true and predicted values of the response simulated for 1000 timesteps between the ground truth and the model predicted values. For these predictions, the testing dataset lies outside of the training dataset both in the case of the autoencoder training and the regression model training. The predictions shown include both the 1st and 10th DOF predictions, in both cases, the prediction performance is initially quite good with performance decreasing with in-creasing simulation time. Whilst the performance does degrade over time, the predictions still remain stable and of a similar magnitude and frequency as the true values. The mean squared error values for each of the predictions shown below were 0.025 and 0.048 respectively. These two examples are indicative of the performance of the method across all 20 degrees of freedom of the system.

\begin{figure}[h!]
\centering
\begin{subfigure}[b]{1\textwidth}
\centering
   \includegraphics[width=140mm]{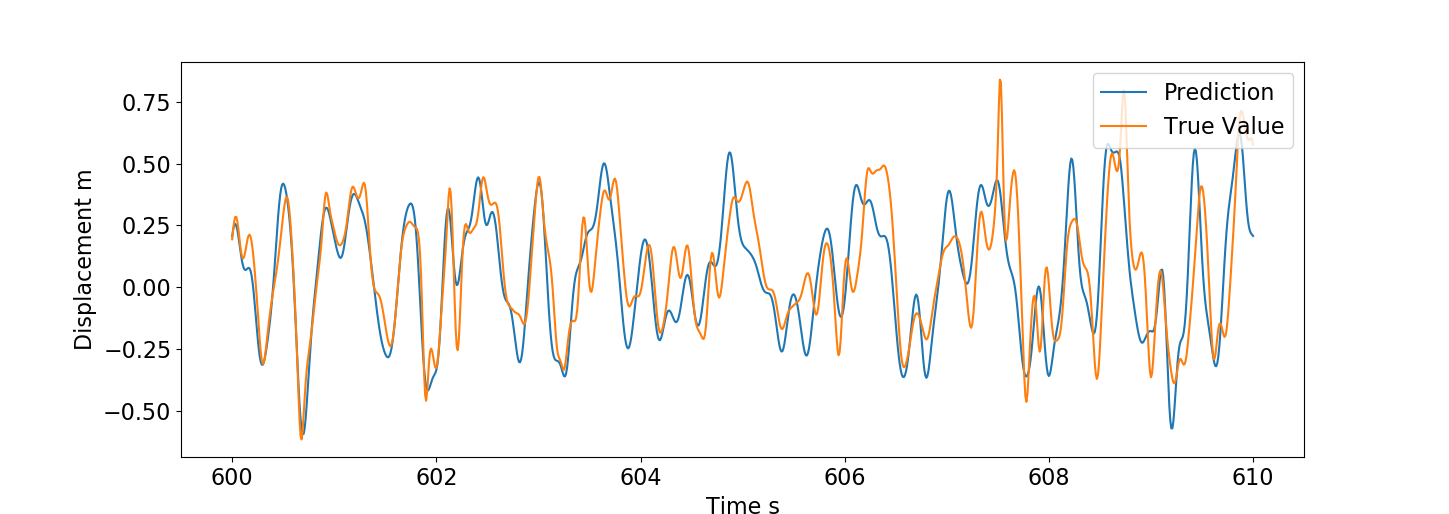}
   \caption{}
   \label{fig:Ng1} 
\end{subfigure}

\begin{subfigure}[b]{1\textwidth}
\centering
   \includegraphics[width=140mm]{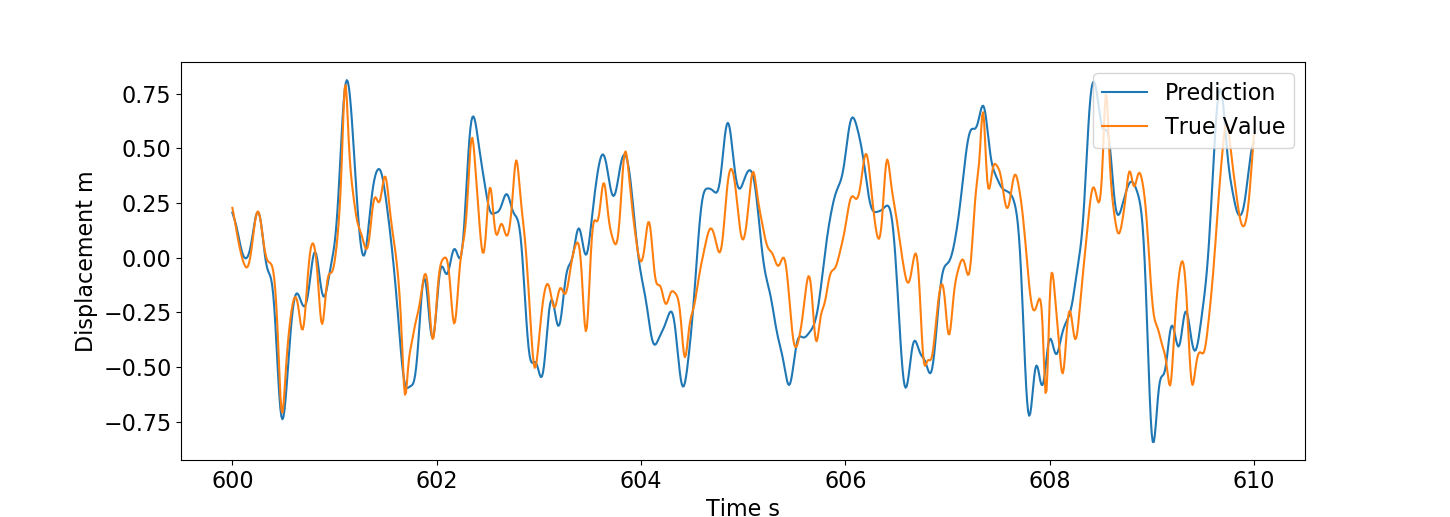}
   \caption{}
   \label{fig:Ng2}
\end{subfigure}
\caption{The final response prediction for both the 1st (a) and 10th (b) DOF for a new force history input}
\end{figure}

\section{Conclusions}
A novel reduced order modelling method was demonstrated, which relies exclusively on input and output data from nonlinear structural systems to allow for time series predictions under arbitrary dynamic forcing. The method makes use of nonlinear normal modes as an efficient reduction basis for nonlinear dynamical systems. It was discussed how these NNMs can be extracted using machine learning methods in an invertible manner and a predictive model was created using an LSTM neural network. The modelling procedure was demonstrated on a 20 DOF nonlinear system featuring cubic nonlinearities. Despite shortcomings in the reconstruction accuracy, which is understandable for this high order of nonlinearity, the proposed method shows considerable promise.

The demonstrated method does offer numerous advantages. It does not assume any a priori form of the nonlinearity considered and hence may be generally applicable to multiple systems. Secondly, due to the nonlinear nature of the NNMs, they should, if sufficient training data is used, be generally valid over the whole state trajectory of a system, as opposed to previous methods based on linear projections.

Future development will concentrate on improving the performance of the regression problem in order to increase the accuracy of the model predictions. In addition, it will be investigated how the amplitude range of the forcing affects the modelling technique and to whether the model can generalise well to broad ranges of input amplitude. Further to this, the method will be tested on larger scale problems, such as finite element simulations  of nonlinear elements or jointed systems. It is also thought that the method may demonstrate greater advantages when applied to larger models, due to a greater number of redundant degrees of freedom. Finally, it is worth noting that variational autoencoders \cite{Kingma2014} can be considered to enforce statistical independence of the extracted components, hence improving predictive capability.

\section*{Acknowledgements}
This work was carried out as part of the ITN project DyVirt and has received funding from the European Union’s Horizon 2020 research and innovation programme under the Marie Skłodowska-Curie grant agreement No 764547

\section*{References}
\begingroup
\renewcommand{\section}[2]{}%
\bibliographystyle{unsrt}
\bibliography{bibmac}
\endgroup

\end{document}